\documentclass{article}
\usepackage{latexsym}
\usepackage{amsthm}
\usepackage{amssymb}
\usepackage{graphicx}

\usepackage{xcolor}
\usepackage{soul}

\newcommand{\Z}{\mathbb{Z}}

\newcommand{\R}{\mathbb{R}}

\newcommand{\x}{{\bf x}}
\newcommand{\bu}{{\bf u}}
\newcommand{\D}{{\bf D}}

\def\Vec{\mathop{\rm Vec}\nolimits}
\newtheorem{theorem}{Theorem}

\newtheorem{lemma}{Lemma}

\title{Some applications of parabolic Dirac operators to the instationary Navier-Stokes problem on conformally flat cylinders and tori in $\R^3$}
\author{
P. Cerejeiras \thanks{Departamento de Matem\'atica, Universidade de Aveiro, P 3810-193 Aveiro, Portugal. E-Mail: {\tt pceres@ua.pt}}
\and
U. K\"ahler \thanks{Departamento de Matem\'atica, Universidade de Aveiro, P 3810-193 Aveiro, Portugal. E-Mail: {\tt ukaehler@ua.pt}}
\and
R.S.~Krau{\ss}har \thanks{Lehrgebiet f\"ur Mathematik und ihre Didaktik, Erziehungswissenschaftliche Fakult\"at, Universt\"at Erfurt, Nordh\"auser Str. 63, D-99089 Erfurt, Germany.  E-mail: {\tt soeren.krausshar@uni-erfurt.de}}}
\begin{document}
\maketitle
\begin{abstract}
In this paper we give a survey on how to apply recent techniques of Clifford analysis over conformally flat manifolds to deal with instationary flow problems on cylinders and tori. Solutions are represented in terms of integral operators involving explicit expressions for the Cauchy kernel that are associated to the parabolic Dirac operators acting on spinor sections of these manifolds.   
\end{abstract}



\section{Introduction}

The  treatment of Navier-Stokes systems is a principal topic in mathematical physics, as they are a main model for describing physical phenomena linked to Newtonian flow from the water flow in a pipe to air flow around a wing. Thus, these systems form the basis of fluid dynamics. The literature addressing these kind of problems is large and abundant ranging from the development of efficient numerical methods up to   theoretical  analysis in the hope to gain some more insight into the structure and nature of  solutions in some special cases. The proof of existence of strong solutions over all times is still open and belongs to the millennium prize problems.
This shows the importance of the field as well as the growing need for further research of particular theoretical aspects of said systems related to fluid dynamics.

\par\medskip\par
In this paper we revisit the three dimensional instationary Navier-Stokes equations  (NSE)  for incompressible fluids
\begin{eqnarray}
-\Delta {\bf u} + \frac{\partial {\bf u}}{\partial t} + ({\bf u} \;{\rm grad})\; {\bf u} + {\rm grad}\; p  &=& {\bf f}, \quad {\rm in}~G\\
{\rm div}\; {\bf u} &=& 0, \quad {\rm in}~G\\
{\bf u}&=&{\bf 0}, \quad {\rm on}~\partial G.
\end{eqnarray}
Here, ${\bf u}$ represents the velocity of the flow, $p$ stands for the pressure and ${\bf f}$ for the specific body force.

\par\medskip\par

Since the 1980s quaternionic analysis, and more generally Clifford analysis, proves to be a powerful tool for the study of this type of non-linear PDE system. In the Stokes problem, which can be regarded as a much simpler stationary and linear version of the NSE, the Laplacian $\Delta = \sum_{i=1}^3 \frac{\partial^2}{\partial x_i^2}$  can be factorized by a linear first order operator, namely the Euclidean Dirac operator $D : = \sum_{i=1}^3 \frac{\partial}{\partial x_i} e_i$ viz $D^2=-\Delta$.
Analogously to the Cauchy-Riemann operator in two dimensions, the Dirac operator is the basis for a rich function theory in higher dimensions. One obtains Cauchy's integral formula in complete analogy with the standard complex Cauchy formula as the basis for the development of further interesting  results, see for instance \cite{DSS} and elsewhere. Similar as in the complex case the resulting function theory provides a refinement of classical harmonic analysis. Furthermore, since the Dirac operator acts on spinor-valued functions, it encodes much more geometrical information than the classical Laplace operator which acts on scalar-valued functions.

\par\medskip\par

The quaternionic operator calculus for elliptic boundary value problems was developed in several works of K. G\"urlebeck, W. Spr\"o{\ss}ig, M. Shapiro, V.V. Kravchenko and many others, see for instance \cite{GS1,GS2,K,KS}. In these works a remarkable number of stationary linear and non-linear boundary value problems have been addressed successfully by means of singular integral operators. The quaternionic calculus actually leads to further new explicit criteria for regularity, existence and uniqueness of the solutions of these systems. Based on those new theoretical results also new numerical algorithms based on discrete version of the quaternionic calculus were developed, see for instance \cite{FGHK,GH}, or the book of M. Mitrea \cite{Mi} in which the study of singular integral operators and Clifford wavelets has successfully been applied to boundary value probems over Lipschitz domains, see also his follow-up work \cite{Mi2}. Also fully analytic representation formulas for the solutions for the Navier-Stokes equations and for Maxwell and Helmholtz systems could be established for some special but important classes of domains \cite{ConKra5,ddsmmas}. A significant advantage of the quaternionic calculus is that its formulae hold universally for all bounded Lipschitz domains, independently of its particular geometry. Furthermore, one gets very convenient analytic representation formulae as well as rather explicit existence and uniqueness criteria. The application of the quaternionic calculus also leads to explicit expressions for the Lipschitz contraction constant for the fixed-point method solving non-linear problems. Based on the explicit knowledge of the contraction constant one obtains useful a-priori and a posteriori estimates on the iterative approximations.

As shown by Sijue Wu in \cite{SW}, the quaternionic analysis calculus turned out be a key ingredient in solving fundamental problems related to the well-posedness of the full 3D water wave problem in Sobolev spaces where the application of well established methods from harmonic and numerical analysis methods did not lead to any success.

About ten years ago as shown in several papers by P. Cerejeiras, U. K\"ahler, F. Sommen, and others cf. e.g. \cite{CK1,CK2} these methods were adapted for dealing with the instationary counterparts of these problems in unbounded Lipschitz domains by means of introducing a parabolic Dirac operator which factorizes the heat operator.
To treat the time-dependent case one adds to the canonical basis elements $e_1,e_2,e_3$   two further basis elements $\mathfrak{f}$ and $\mathfrak{f}^{\dagger}$  which satisfy $\mathfrak{f}^2 = (\mathfrak{f}^{\dagger})^2 = 0$. The additional elements are often called a Witt basis and they allow for the inclusion of the time dimension. 

By means of a general positive real $k > 0$ we consider the (modified) parabolic Dirac operator
$$
D_{{\bf x},t;k}^{\pm} := \sum\limits_{j=1}^3 e_j \frac{\partial }{\partial x_j} + \mathfrak{f} \frac{\partial }{\partial t} \pm k^2 \mathfrak{f}^{\dagger}.
$$
This operator factorizes the generalized heat operator viz $(D_{{\bf x},t;k}^{\pm})^2 = -\Delta \pm k^2 \frac{\partial }{\partial t}$.
Similarly to the elliptic case, for this operator one can also introduce  adequate analogues of the Teodorescu transform, the regular and singular Cauchy transform and the Bergman projection operator, as proposed in \cite{Be,CK2,CV2009,KV}. This adapted operator calculus allows us to treat the time-dependent versions of the PDEs studied earlier over time varying domains in a rather elegant way.

\par\medskip\par

A further recent line of investigation consists in developing possible extensions of this operator calculus to handle such PDE on curved spaces and manifolds. If we want to study for instance weather cast problems, then one appropriate model consists of working with the Navier-Stokes system on the sphere. The latter then involves spherical versions of the Laplacian and the Dirac operator, as proposed for instance by W. Spr\"o{\ss}ig in \cite{Sprweather}. In comparison with the earth radius the atmosphere has a negligible small thickness, so that one deals at  first glance with a flow problem on a sphere. To apply the quaternionic operator calculus to the spherical case, one has to adapt the integral operators in a geometric appropriate way, namely the Euclidean Cauchy kernel has to be substituted by its spherical analogue. However, the representation of the solution again consists of the same types of integral operators as introduced in the Euclidean case. One simply has to compute the kernel functions for the new geometry and to replace the Euclidean kernels by the properly  adapted versions of these kernels.

\par\medskip\par

In this paper we want to outline how one can adapt the quaternionic operator calculus to study time-dependent Navier-Stokes problems in the more general context of conformally flat spin manifolds that arise by factorizing out some simply connected domain by a discrete Kleinian group. Here, we treat conformally flat spin cylinders an tori as an illustrative example. This underscores the universality of our approach. Furthermore, our approach has the advantage that the representation formulas and criteria can directly be generalized to the $n$-dimensional setting, just by replacing the quaternionic operators by their corresponding Clifford algebra valued ones, such as suggested in \cite{CK2} for the Navier-Stokes system, in which all formulas remain with the same structure. We are able to construct the parabolic Cauchy kernel on these manifolds explicitly. This consequently opens the door to apply the iterative computation algorithm to compute the solutions.

\section{Preliminaries}

\subsection{Quaternionic function theory}

Let $\{e_1,e_2,e_3\}$ denote the standard basis of the Euclidean vector space $\mathbb{R}^3$. To endow the space $\mathbb{R}^3$
with an additional multiplicative structure, we embed it into the Hamiltonian algebra of real quaternions, denoted by $\mathbb{H}$. A
quaternion is an element of the form $$x=x_0+{\bf x} := x_0 + x_1 e_1 + x_2 e_2 + x_3 e_3$$ where $x_0,\ldots,x_3$ are real numbers.
$x_0$ is called the real part of the quaternion $x$ and will be denoted by $\Re(x)$ while ${\bf x},$ or $\Vec(x),$ denotes the vector part of $x$. In the quaternionic setting the standard unit vectors play the role of imaginary units, i.e., we have 
$$e_1e_2e_3 = e_i^2 = -1, \quad i=1,2,3.$$
The generalized anti-automorphism \textit{conjugation} in $\mathbb{H}$ is defined by
$$\overline{1} = 1, \quad \overline{e_i} = - e_i, ~i=1,2,3,\quad \overline{ab} = \overline{b} \; \overline{a}.$$ The Euclidean norm in $\R^4$ induces a norm on the whole quaternionic algebra as $|a|:= \sqrt{\sum_{i=0}^3
a_i^2}$.\\

In what follows let $G \subset \R^3$ be a bounded domain with a sufficiently smooth boundary $\Gamma=\partial G.$ A quaternionic function $f : G \subset \R^3 \rightarrow \mathbb{H}$ has a representation $$\x \mapsto f(\x) = f_0(\x) + {\bf f}(\x) := f_0(\x)+  \sum_{i=1}^3 f_i(\x) e_i,$$ with $\R-$valued components $f_i.$ Properties like continuity, etc. are understood coordinatewisely. Now, the additional multiplicative structure of the quaternions allows to describe all $C^1$-functions $ {\bf f}:\mathbb{R}^3 \rightarrow
\mathbb{R}^3$ that satisfy both ${\rm div}\;{\bf f} = 0$ and ${\rm rot}\;{\bf f}=0$ in a compact form as null-solutions of the three-dimensional Dirac operator $${\bf D}:=\sum_{i=1}^3 \frac{\partial }{\partial x_i}
e_i.$$ 
This operator is nothing else than the Atiyah-Singer-Dirac operator that arises in a natural way from the Levi-Civita connection in the context of Riemannian spin manifolds. In the Euclidean $3D$-space it simplifies to the form above. More important, the Euclidean Dirac operator coincides with the usual gradient operator when applied to a scalar-valued function. This motivates the study of \textit{monogenic functions.} A real differentiable function $f: G \subset \R^3 \rightarrow \mathbb{H}$ is called left quaternionic monogenic, or left quaternionic holomorphic, in $G$ if one has $\D f=0$ in $G.$ Since the Euclidean Dirac operator factorize up to signal the (3D) Laplace operator, that is $\Delta f = -\D^2 f,$ we have that  every real component of a left monogenic function is again a harmonic function. Conversely, following e.g. \cite{DSS}, if $f \in C^2$ is a solution of the Laplace operator in $G$, then in any open ball $B(\tilde{\x},r) \subset G$ there exist two left monogenic functions, $f_0$ and $f_1$, such that $f = f_0 + {\bf x} f_1$ holds in $B(\tilde{\x},r)$. This property allows to treat harmonic functions in terms of null solutions of $\D.$ It is also the starting point for the construction of analogues to several well known theorems of complex analysis. For more details on quaternionic functions and operator theory, we refer the reader for instance to \cite{DSS,GS1,GS2}.

\subsection{The instationary case}

To treat time dependent problems in $\R^3$ we follow the ideas of \cite{CK2}. First, we introduce two additional basis elements $\mathfrak{f}$ and $\mathfrak{f}^{\dagger}$ satisfying to
\begin{equation} 
\mathfrak{f} \mathfrak{f}^{\dagger} + \mathfrak{f}^{\dagger} \mathfrak{f} = 1,\quad \mathfrak{f}^2 = (\mathfrak{f}^{\dagger})^2 = 0, \label{Eq:00a}
\end{equation}
and which interact with the existent elements of the basis of $\R^3$ as
\begin{equation} \mathfrak{f} e_j = e_j \mathfrak{f} =0, \quad \mathfrak{f}^{\dagger} e_j = e_j \mathfrak{f}^{\dagger} = 0.
\label{Eq:00b}
\end{equation}
We construct the (dual) parabolic Dirac operators given by
$$
D_{{\bf x},t}^{+} := \D+ \mathfrak{f} \frac{\partial }{\partial t} + \mathfrak{f}^{\dagger}, \quad 
D_{{\bf x},t}^{-} := \D+ \mathfrak{f} \frac{\partial }{\partial t} - \mathfrak{f}^{\dagger}.$$ 
We remark that, based on (\ref{Eq:00a}) and (\ref{Eq:00b}) these operators satisfy $(D_{{\bf x},t}^{\pm})^2 = -\Delta \pm \frac{\partial }{\partial t},$ that is to say, they factorize the \textit{heat operator}. Null solutions of the parabolic Dirac operator $D_{{\bf x},t}^{+}$ are called (left) parabolic monogenics  (resp. dual parabolic monogenics if solutions of $D_{{\bf x},t}^{-} f = 0$). \\

Suppose now that $G$ is a space-time varying bounded Lipschitz domain $G \subset \R^3 \times \R^+.$ In what follows we define $W_2^{s,l}(G)$ as the parabolic Sobolev space of $L_2(G)$ where $s$ is the regularity parameter with respect to ${\bf x}$ and $l$ the regularity parameter with respect to $t$. Using the Stokes theorem we get 
$$\int_G \Big[  \Big(g(\D +\mathfrak{f} \partial_t)\Big) f + g \Big((\D + \mathfrak{f} \partial_t)   \Big) f \Big] d\x = \int_{\Gamma} g d\sigma_{\x,t} f, $$ where $d\sigma_{{\bf x},t} = \Big( \D + \mathfrak{f} \frac{\partial }{\partial t} \Big) \rfloor d\x dt$ is the contraction of the homogenous operator $  \D + \mathfrak{f} \partial_t$ with the volume element $d\x dt$. Hence, this leads to the Stokes integral formula involving out parabolic Dirac operators, namely
\begin{equation} \label{Eq:00c}
\int_G \Big[  (g D_{{\bf x},t}^{-}) f + g (D_{{\bf x},t}^{+}f ) \Big] d\x = \int_{\Gamma} g d\sigma_{\x,t} f.
\end{equation}

Moreover, the fundamental solution to the dual operator $D_{{\bf x},t}^{-}$ has the form 
$$
E_-({\bf x},t) = \frac{H(t) \exp(-\frac{|{\bf x}|^2}{4t})}
{(2 \sqrt{\pi t})^3}
\Big(\frac{\x}{2t}  + \mathfrak{f}(\frac{3}{2t}  + \frac{|{\bf x}|^2}{4t^2})- \mathfrak{f}^{\dagger}\Big),
$$
where $H(\cdot)$ stands for the usual Heavyside function. Replacing the fundamental solution in (\ref{Eq:00c}) we obtain the Borel-Pompeiu integral formula. 

\begin{theorem}(see \cite{CK2,CV2009})\label{bp} 
Let $G \subset \mathbb{R}^3 \times \R^+$ be a bounded Lipschitz domain with a strongly Lipschitz boundary $\Gamma=\partial G$. 

Then for all $u \in W_{2}^{1,1}(G)$ we have
$$
\int_{\Gamma} E({\bf x}-{\bf y},t-\tau) d\sigma_{{\bf x},t} u({\bf x},t) = u({\bf y},\tau) + \int_{G} E({\bf x}-{\bf y},t-\tau) (D_{{\bf x},t}^{+}u)({\bf x},t) d\x dt.
$$ 
\end{theorem}

Whenever $u \in $ Ker $D_{{\bf x,t}}^{+}$ one obtains the following version of Cauchy's integral formula for parabolic monogenic functions $$
u({\bf y},\tau) = \int_{\Gamma} E({\bf x}-{\bf y},t-\tau) d\sigma_{{\bf x},t} u({\bf x},t).
$$
Again, following the above cited works, we can introduce the parabolic Teodorescu transform and the Cauchy transform by 
\begin{eqnarray*}
T_G u({\bf y},\tau) &=& \int_G E({\bf x}-{\bf y},t-\tau) u({\bf x},t) d\x  dt,\\
F_{\Gamma} u({\bf y}, \tau) &=& \int_{\Gamma} E({\bf x}-{\bf y},t-\tau) d\sigma_{{\bf x},t} u({\bf x},t). 
\end{eqnarray*} On the one hand we have $D_{{\bf x},t}^+T_G u = u,$ that is, the parabolic Teodorescu operator is the right inverse of the parabolic Dirac operator. On the other hand, and analogously to the Euclidean case we can rewrite the Borel-Pompeiu formula in the form 

\begin{lemma}\label{Lem:01} Let $u \in W_{2}^{1,0}(G)$. Then $$T_G D_{{\bf x},t}^+u = u - F_{\Gamma} u.$$   
\end{lemma} 

The space $L_2(G)$ can be decomposed into the direct sum of the subspace of parabolic monogenics in $G$ and its
complement. 

\begin{theorem} (Hodge decomposition). Let $G \subseteq \mathbb{R}^3 \times \R^+$ be a bounded Lipschitz domain.
Then $$L_2(G) = \Big( L_2(G) \cap {\rm Ker}D_{{\bf x},t}^+ \Big) \oplus D_{{\bf x},t}^+ {\stackrel{\circ}{W}}_{2}^{1,1}(G).$$
where $L_2(G) \cap {\rm Ker} D_{{\bf x},t}^+ =: B(G)$ is the Bergman space of parabolic monogenic functions, and where ${\stackrel{\circ}{W}}_{2}^{1,1}(G)$ is the subspace of all $f \in W_{2}^{1,1}(G)$ with vanishing boundary data. 
\end{theorem} 

Proofs of the above results can be found in \cite{CK2,CV2009}. \\

{\bf Remark}: Due to the exponential decrease of the fundamental solution, the operator $T_G$ remains a $L^2(G)$ bounded operator also if $G \subset \R^3 \times \R^+$ is unbounded. The application of the add-on term as proposed in \cite{CK1} for the Teodorescu transform associated to the usual spatial Dirac operator ${\bf D}$ is not necessary in the parabolic setting. \\

For our purpose we need the more general parabolic Dirac operator, used for instance in \cite{Be,CV2009,KV}, having the form 
$$
D_{{\bf x},t,k}^{\pm} := \D + \mathfrak{f} \frac{\partial }{\partial t} \pm k \mathfrak{f}^{\dagger} =\sum\limits_{j=1}^3 e_j \frac{\partial }{\partial x_j} + \mathfrak{f} \frac{\partial }{\partial t} \pm k \mathfrak{f}^{\dagger}
$$ 
for a positive real $k \in \R$. This operator factorizes the second order operator 
$$
(D_{{\bf x},t,k}^{\pm})^2 = -\Delta \pm k^2 \frac{\partial }{\partial t}
$$
and has very similar properties as the previously ones. Their null-solutions are called parabolic $k$-monogenic (resp. dual parabolic $k$-monogenic) functions. 

Adapting from \cite{Be,CV2009}, the fundamental solution to $D_{{\bf x},t,k}^{+}$ turns out to have the form 
$$
E({\bf x},t;k) = \sqrt{k}\frac{H(t) \exp(-\frac{k|{\bf x}|^2}{4t})}
{(2 \sqrt{\pi t})^3}
\Big(\frac{k}{2t} \sum\limits_{j=1}^3 e_j x_j + \mathfrak{f}(\frac{3}{2t} 
+ \frac{k|{\bf x}|^2}{4t^2})+k^2 \mathfrak{f}^{\dagger}\Big).
$$

\par\medskip\par
In what follows 
${\bf P} : L_2(G) \rightarrow B(G):= L_2(G) \cap {\rm Ker} D_{{\bf x},t;k}^+$ denotes the
orthogonal Bergman projection while ${\bf Q} : L_2(G) \rightarrow
D_{{\bf x},t}^+ {\stackrel{\circ}{W}}_{2}^{1,1} (G)$ stands for the projection
into the complementary space in all that follows. One has ${\bf Q}
= {\bf I} - {\bf P}$. Here ${\bf I}$ stands for the identity operator.\\

The Bergman space of parabolic $k$-monogenic functions is a Hilbert space with a uniquely defined reproducing kernel function, the so-called the parabolic $k$-monogenic Bergman kernel denoted by $B({\bf x},{\bf y};t,\tau)$. The
orthogonal Bergman projection ${\bf P}: L_2(G) \rightarrow
B(G)$ is given by the convolution with the Bergman
kernel
$$
({\bf P} u)(\x, t) = \int_G B({\bf x},{\bf y};t,\tau) u({\bf y},\tau) d{\bf y}d\tau, \quad\quad u \in L_2(G).
$$
In particular, one has $({\bf P} u)({\bf x},t) = u({\bf x},t)$ for all $u \in B(G)$.

\section{The Navier-Stokes equations shortly revisited in quaternions}

In the classical vector analysis calculus the in-stationary Navier-Stokes equations have the form (again, we assume here viscosity $\nu = 1$)
\begin{eqnarray}
-\Delta {\bf u}  + \frac{\partial {\bf u}}{\partial t} +  ({\bf u} \;{\rm grad})\; {\bf u} + {\rm grad}\; p  &=& {\bf f}, \quad {\rm in}~G\\
{\rm div}\; {\bf u} &=& 0, \quad {\rm in}~G\\
{\bf u}&=&{\bf 0}, \quad {\rm on} ~\partial G
\end{eqnarray}
To apply the quaternionic integral operator calculus to solve these equations one first expresses this system in the quaternionic language, as done in \cite{CK2}. \\

First we recall that for a time independent quaternionic function of type 
$${\bf f} : \R^3 \to \R^3,\quad  \mbox{with }\x \mapsto {\bf f}(\x)$$ we have  
$$
\D {\bf f} = {\rm rot}\;{\bf f} - {\rm div}\; {\bf f}.$$
 Hence, the divergence of a vector field ${\bf f}$ can be expressed as div ${\bf f} = \Re({\bf D}{\bf f}$). 


In a similar way, for a scalar valued function $p : \R^3 \to \R, ~\x \mapsto p(\x),$  we have $${\bf D} p = {\rm grad}\;p.$$

Finally, we recall that the three dimensional Euclidean Laplacian $\Delta=\sum_{i=1}^3 \frac{\partial^2}{\partial x_i^2}$ can be expressed in terms of the Dirac operator as $\Delta=-{\bf D}^2.$
\par\medskip\par
Next, we assume that the vector-field is time-dependent, that is, ${\bf u} = {\bf u}(\x, t).$ Applying the formulas from the preceding section we can express the heat operator $- \Delta {\bf u}  + \frac{\partial {\bf u}}{\partial t}$ in the form 
\begin{eqnarray*}
- \Delta {\bf u}  + \frac{\partial {\bf u}}{\partial t} &=& (D_{{\bf x},t}^+)^2 {\bf u}\\
\end{eqnarray*}
Thus, the original system for a time-dependent vector-field ${\bf u} = {\bf u}(\x, t)$ can be reformulated in the following way: 
\begin{eqnarray}
(D_{{\bf x},t}^+)^2 {\bf u} + \Re({\bf u} \;{\bf D})\; {\bf u} + {\bf D}\; p  &\!\!\!=\!\!& {\bf f} \;\;{\rm in}\;G\\
\Re({\bf D} {\bf u}) &\!\!\!=\!\!& 0 \;\;{\rm in}\;G\\
{\bf u}&=&{\bf 0},\; {\rm at}\;\partial G.
\end{eqnarray}
The strategy for the resolution of this system is to apply the previously introduced hypercomplex integral operators in order to get iterative  formulas for the velocity ${\bf u}$ and the pressure $p$. 

\section{The linear case}

In this section we briefly recall how the quaternionic calculus can be applied to set up analytic solutions for the special case in which the convective term $({\bf u} \;{\rm grad})\; {\bf u}$ is negligibly small. Hence, we assume $p$ and the external source ${\bf f}$ in $L_2(G)$ and ${\bf u} \in {\stackrel{\circ}{W}}_{2}^{1,1}(G).$ 

~

Under these assumptions the instationary viscous Navier-Stokes equations take the simplified form 
\begin{eqnarray}
(D_{{\bf x},t}^+)^2 {\bf u}  + {\bf D}\; p  &=& {\bf f} \;\;{\rm in}\;G  \label{v}  \\
\Re({\bf D} {\bf u}) &=& 0 \;\;{\rm in}\;G  \label{v2} \\
{\bf u}&=&{\bf 0}\;{\rm at}\;\partial G. 
\end{eqnarray} 

The velocity ${\bf u}$ and the pressure $p$ can now be computed from this system using the (modified) parabolic Teodorescu operator.  
Applying the parabolic Teodorescu operator to (\ref{v}) leads to the equation 
\begin{equation}\label{eq1}
(T_G D_{{\bf x},t}^+)(D_{{\bf x},t}^+ {\bf u}) + T_G {\bf D} p = T_G ({\bf f}).
\end{equation}
Now, we apply Lemma \ref{Lem:01} (Borel-Pompeiu formula) to~(\ref{eq1}).  This leads to 
\begin{equation}\label{eq2}
(D_{{\bf x},t}^+{\bf u} - F_{\Gamma} D_{{\bf x},t}^+ {\bf u}) + T_G {\bf D} p = T_G ({\bf f}).
\end{equation}
Using the orthogonal Bergman projector ${\bf Q}$ yields 
\begin{equation}\label{eq3}
({\bf Q}D_{{\bf x},t}^+{\bf u} - {\bf Q} F_{\Gamma} D_{{\bf x},t}^+ {\bf u}) + {\bf Q} T_G {\bf D}p = {\bf Q} T_G ({\bf f}).
\end{equation}
Since the Cauchy integral operator maps $L_2(G)$  onto  $L_2(G) \cap $Ker $D_{{\bf x},t}^+$ and ${\bf u} \in {\stackrel{\circ}{W}}_{2}^{1,1}(G)$ we get that $F_\Gamma D_{{\bf x},t}^+ {\bf u}$ is a left parabolic monogenic function that is, ${\bf Q} F_{\Gamma} D_{{\bf x},t}^+ {\bf u} = 0$.

Therefore, equation~(\ref{eq3}) simplifies to 
\begin{equation}\label{eq4}
{\bf Q}D_{{\bf x},t}^+{\bf u} + {\bf Q}T_G  {\bf D} p = {\bf Q}T_G ({\bf f}).
\end{equation} At this point we remark that $T_G$ is the right inverse to $D_{{\bf x},t}^+$ but not to ${\bf D}$! 

~

We apply again the Teodorescu transform to equation~(\ref{eq4}):
\begin{equation}\label{eq5}
T_G {\bf Q} D_{{\bf x},t}^+{\bf u} + T_G{\bf Q}T_G {\bf D}p = T_G{\bf Q}T_G ({\bf f}).
\end{equation}
First, we observe that $T_G {\bf Q} D_{{\bf x},t}^+ {\bf u} = T_G D_{{\bf x},t}^+ {\bf u}$, because $D_{{\bf x},t}^+ {\bf u} \in \;im({\bf Q})$. 
Applying again Lemma \ref{Lem:01} leads to 
\begin{equation}\label{eq6}
 {\bf u} -   F_{\Gamma} {\bf u} + T_G{\bf Q} T_G {\bf D} p = T_G{\bf Q}T_G ({\bf f}).
\end{equation}
Since  ${\bf u}|_{\Gamma} = {\bf 0}$, we get that  $F_{\Gamma} {\bf u}$  vanishes and we do obtain the following representation formula for the velocity field ${\bf u}$:
\begin{equation}\label{velocity}
{\bf u} = T_G {\bf Q} T_G \left( {\bf f} -  {\bf D} p \right).
\end{equation}
The pressure $p$ can be obtained from the continuity equation (\ref{v2}). Indeed, inserting the solution ${\bf u}$ obtained in (\ref{velocity}) into (\ref{v2}) leads to 
\begin{equation}\label{pressure}
\Re({\bf Q}T_G {\bf D} p)= \Re ({\bf Q} T_G {\bf f}).
\end{equation}  
Thus, we have obtained the following representation formulas for the solutions of the instationary viscous Navier-Stokes equations in the case of a negligibly small convective term:
\begin{theorem}\label{th3} (Representation theorem). Suppose that $p \in L_2(G), {\bf u} \in {\stackrel{\circ}{W}}_{2}^{1,1}(G)$ Then the solutions can be represented in the form 
\begin{eqnarray}
\Re({\bf Q}T_G{\bf D} p) & = & \Re ({\bf Q} T_G  {\bf f}) \label{LP}\\
{\bf u} & = &  T_G {\bf Q} T_G {\bf f} - T_G {\bf Q} T_G {\bf D} p \label{LV},
\end{eqnarray}
\end{theorem}
The pressure is uniquely determined from (\ref{LP}) up to a constant. Given the solution for the pressure $p$, (\ref{LV}) gives the solution for ${\bf u}.$ Hence, the original system is solvable by application of the integral operators $T_G$ and ${\bf Q}$. 
\par\medskip\par
Both the Teodorescu and Cauchy integral operators have a universal integral kernel for all bounded domains, namely the Cauchy kernel $E = E({\bf x}-{\bf y};t -\tau)$. Also, the Bergman projectors can be expressed by the algebraic relation
$$
{\bf P} = F_{\Gamma}(tr_{\Gamma} T_G F_{\Gamma})^{-1} tr_{\Gamma} T_G,
$$
where $tr_{\Gamma}$ is the usual trace operator, or the restriction to the boundary of the domain. Furthermore, the above scheme is extendable to the case of $k-$monogenics, that is, when (\ref{v}) is given as $(\D_{{\bf x}, t; k})^2 {\bf u}+\D p = {\bf f}.$ See \cite{GS1, CK2} for more details. 

\section{The case of a non negligenciable convective term} 

Now, we turn our attention to the more complicated case in which the flow is still viscous ($\nu =1$) but the non-linear convective term $({\bf u}~{\rm grad}){\bf u} = \Re({\bf u}{\bf D}){\bf u}$ is no longer negligibly small.  
\par\medskip\par
First, we observe that the reasoning and arguments used in Section~4 are still valid. Hence, we get the following equations for the velocity ${\bf u}$ and pressure $p$, that is, 
\begin{equation}\label{velocitynonlinear}
{\bf u} =  T_G {\bf Q} T_G \Big[ {\bf f} -\Re({\bf u}{\bf D}){\bf u}\Big] - T_G {\bf Q} T_G {\bf D} p,
\end{equation}
and, from inserting this solution into the continuity equation (\ref{v2}), we obtain  
\begin{equation}\label{pressurenonlinear}
\Re({\bf Q}T_G {\bf D} p)= \Re\Big[{\bf Q} T_G \left( {\bf f} - \Re({\bf u}{\bf D}){\bf u} \right) \Big].
\end{equation}

Now, we apply the following fixed point algorithm in order to iteratively compute both solution ${\bf u}$ and pressure $p,$ departing from an arbitrary ${\bf u}_0$ (for the time being, no conditions will be imposed here):
\begin{eqnarray*}
\Re({\bf Q}T_G {\bf D} p_n)&=&
\Re\Big[{\bf Q} T_G\Big(  {\bf f} - \Re({\bf u}_{n-1}{\bf D}){\bf u}_{n-1} \Big) \Big],\\
{\bf u}_n &=&  T_G {\bf Q} T_G \Big[{\bf f}-\Re({\bf u}_{n-1}{\bf D}){\bf u}_{n-1}\Big] - T_G {\bf Q} T_G {\bf D} p_n
\end{eqnarray*} for $n=1, 2, \ldots$
  
\par\medskip\par
The following lemma (c.f. \cite{CK2}) gives the conditions under which the above proposed fixed point algorithm does converge to a unique solution $({\bf u}, p):$
\begin{lemma}\label{helpl} 
Suppose that ${\bf u} \in {\stackrel{\circ}{W}}_2^{s,l}(G) \cap Ker D_{{\bf x},t},$ with $s,l \ge 1,$ and $p \in L_2(G)$ are pairwise solutions of  (\ref{velocitynonlinear}) and (\ref{pressurenonlinear}). Then, the following  estimate holds:
\begin{equation}
\|D_{{\bf x},t}^+ {\bf u}\|_2  + \| {\bf Q} p\|_2  \le \sqrt 2 \| T_G M({\bf u})\|_2,
\end{equation}    
where 
$M({\bf u}) := \Re({\bf u}{\bf D}){\bf u} - {\bf f}$ and $\| \cdot \|_2$ stands for the $L_2-$norm.
\end{lemma}

In fact, for $p \in W^{1,1}_2(G)$ we have 
$$T_G {\bf Q} T_G \D p =T_G {\bf Q} (p - F_\Gamma p) =  T_G {\bf Q} p,$$ since  $F_\Gamma p$ is in ${\rm im} {\bf P}.$ Applying $D_{{\bf x},t}^+$ to this equation gives (recall that $D_{{\bf x},t}^+$ is a left inverse for  $T_G$) 
$$D_{{\bf x},t}^+ ( T_G {\bf Q} T_G \D p ) =   {\bf Q} p.$$
Since $W^{1,1}_2(G)$ is dense in $L_2(G)$ this leads to
$D_{{\bf x},t}^+ ( T_G {\bf Q} T_G \D p ) =   {\bf Q} p$ for all $p \in L_2(G),$
and we obtain from (\ref{pressurenonlinear}) 
$$D_{{\bf x},t}^+ {\bf u} =  {\bf Q} T_G M({\bf u}) - {\bf Q} p.$$

By the orthogonality between $D_{{\bf x},t}^+ {\bf u}$ and ${\bf Q} p,$ we have 
$$ \| D_{{\bf x},t}^+ {\bf u}\|_2 + \| {\bf Q} p\|_2 \le \sqrt 2 \| {\bf Q} T_G M({\bf u})\|_2 = \sqrt 2 \| T_G M({\bf u})\|_2.$$

\par\medskip\par

Starting with $\bu_0 \in {\stackrel{\circ}{W}}_{2}^{1,1}(G)$ we generate the iteration pairs $(p_n,{\bf u}_n) \in L_2(G) \times {\stackrel{\circ}{W}}_{2}^{1,1}(G).$ Moreover, by (\ref{pressurenonlinear}) we get
\begin{eqnarray*}
\|{\bf u}_n - {\bf u}_{n-1}\|_{W_{2}^{1,1}}  &\le &  \| T_G {\bf Q} T_G [M({\bf u}_{n-1}) - M({\bf u}_{n-2}) ] \|_{W_{2}^{1,1}} + \| T_G{\bf Q} (p_n - p_{n-1}) \|_{W_{2}^{1,1}} \\
& \le &  2C_1 \| M({\bf u}_{n-1}) - M({\bf u}_{n-2})  \|_{W_{2}^{-1,-1}}, 
\end{eqnarray*}
where $C_1 := \|T_G {\bf Q} T_G\|.$ 

Moreover, according to \cite{CK1} Lemma 4.1 for all $\bu \in {\stackrel{\circ}{W}}_{2}^{1,1}(G)$ there exists a constant $C_2$ such that  
$$
\|\Re({\bf u}{\bf D}){\bf u}\|_{W_{2}^{-1,-1}} \le C_2 \|{\bf u}\|_{W_{2}^{1,1}}^2,
$$ 
so that the previous estimate becomes   
$$
\|{\bf u}_n - {\bf u}_{n-1}\|_{W_{2}^{1,1}} \le 2 C_1 C_2  \Big( \|{\bf u}_{n-1} \|_{W_{2}^{1,1}} + \|{\bf u}_{n-2}\|_{W_{2}^{1,1}}    \Big) \|{\bf u}_{n-1} - {\bf u}_{n-2}\|_{W_{2}^{1,1}}.
$$
Next, we need to prove that the energy of our solutions ${\bf u}_n$ decreases, that is to say, $\|{\bf u}_{n}\|_{W_{2}^{1,1}} \le \|{\bf u}_{n-1}\|_{W_{2}^{1,1}}$. First, we observe that 
\begin{eqnarray*}
\|{\bf u}_{n}\|_{W_{2}^{1,1}} & \le & \| T_G {\bf Q} T_G {\bf u}_{n-1}\|_{W_{2}^{1,1}} + \|T_G {\bf Q} p_n\|_{W_{2}^{1,1}} \\
& \le & 2 C_1 C_2 \|{\bf u}_{n-1}\|_{W_{2}^{1,1}}^2+ 2C_1 \|{\bf f}\|_{2}.
\end{eqnarray*}
Hence, $\|{\bf u}_{n}\|_{W_{2}^{1,1}} \le \|{\bf u}_{n-1}\|_{W_{2}^{1,1}}$ whenever 
$$2 C_1 C_2 \|{\bf u}_{n-1}\|_{W_{2}^{1,1}}^2+ 2C_1 \|{\bf f}\|_{2} \leq \|{\bf u}_{n-1} \|_{W_{2}^{1,1}}$$ which leads to
$$
\Big(\|{\bf u}_{n-1}\|_{W_{2}^{1,1}} - \frac{1}{4C_1C_2}\Big)^2 \le \frac{1}{16 C_1^2 C_2^2} - \frac{1}{C_2}\|{\bf f}\|_2.$$

Now, if $\|{\bf f}\|_2 \le \frac{1}{16C_1^2C_2}$, then the previous in-equation can be written as
$$\frac{1}{4C_1C_2} -W \leq  \|{\bf u}_{n-1}\|_{W_{2}^{1,1}} \leq \frac{1}{4C_1C_2} + W,$$
where $W := \sqrt{\frac{1}{16 C_1^2 C_2^2} - \frac{1}{C_2}\|{\bf f}\|_2}.$ This finally leads to an estimate on the Lipschitz constant 
\begin{eqnarray*}
\|{\bf u}_n - {\bf u}_{n-1}\|_{W_{2}^{1,1}} & \le & 2 C_1 C_2  \Big( \|{\bf u}_{n-1} \|_{W_{2}^{1,1}} + \|{\bf u}_{n-2}\|_{W_{2}^{1,1}}    \Big) \|{\bf u}_{n-1} - {\bf u}_{n-2}\|_{W_{2}^{1,1}} \\
& \le &   \Big( 1+2 C_1 C_2 W  \Big) \|{\bf u}_{n-1} - {\bf u}_{n-2}\|_{W_{2}^{1,1}},
\end{eqnarray*}
of the form
$$
L := 1 - 4 C_1 C_2 W  < 1.
$$ 
Summarizing 
\begin{theorem} (cf. \cite{CK2} p. 1723).\\
The iteration method converges for each starting point ${\bf u}_{0} \in {\stackrel{\circ}{W}}_{2}^{1,1}(G) \cap  ker D^+_{{\bf x},t}$ with
$$
\|{\bf u}_{0}\|_{W_{2}^{1,1}}  \le \min \left( \frac{1}{2C_1C_2}, \frac{1}{4C_1C_2}+W \right)
$$
and $W:= \sqrt{\frac{1}{16 C_1^2 C_2^2} - \frac{1}{C_2}\|{\bf f}\|_2}.$
\end{theorem}

\par\medskip\par

{\bf Remarks}: we get the pressure $p$ up to an additive constant when $G$ is bounded, and we get uniqueness of $p$ when $G$ is unbounded.  Also, and as explained in \cite{BGSS} in the time independent case we can replace the Teodorescu transform by a simpler primitivation operator whose evaluation requires less computational steps. 

\section{The Navier-Stokes equations in the more general context of some conformally flat spin $3$-manifolds}

One further advantage of using the quaternionic operator calculus consists in the fact that the results and representation formulas presented in the previous sections can easily be carried over to the treatment of analogous boundary value problems within the more general context of conformally flat spin manifolds, of which the Euclidean space $\R^3$ is just the simplest example. This is due to the fact that the formulas presented in the previous sections have geometrically a very universal character. 

\par\medskip\par

Recalling for example from the classical paper \cite{Kuiper} a conformally flat $3$-manifold is a Riemannian $3$-manifold that has a vanishing Weyl tensor. In dimensions $n \ge 3$ these are exactly those Riemannian manifolds that have atlasses whose transition functions are M\"obius transformations. 

As also pointed out in \cite{Kuiper} one way of constructing examples of conformally flat manifolds is to factor out a subdomain $U$ of either the sphere $S^{3}$ or $\R^{3}$ by a Kleinian subgroup $\Gamma$ of the M\"{o}bius group where $\Gamma$ acts totally discontinuously on $U$. This gives rise to the conformally flat manifold $U / \Gamma$. In the original paper by N.H.~Kuiper it is shown that the universal cover of a conformally flat manifold admits a development (i.e. a local conformal diffeomorphism) into $S^3$. The class of conformally flat manifolds of the form $U/\Gamma$ are exactly those for which this development is a covering map $\tilde{U} \rightarrow U \subset S^3$. 

\par\medskip\par

Examples of such manifolds are for example $3$-tori, cylinders, real projective space and the hyperbolic manifolds $H^{+}/\Gamma_p[N]$ for an integer $N \ge 2$ where $H^+:=\{{\bf x} \in \R^3 \mid x_3 > 0$\} and where $\Gamma_p[N]$ is a principal congruence arithmetic subgroup of level $N$ of the hypercomplex modular group $\Gamma_p$. The latter one is generated by the Kelvin inversion (i.e. the reflection at the unit sphere ${\bf x} \mapsto -{\bf x}/|{\bf x}|^2$) and by the translation operations ${\bf x} \mapsto {\bf x} + e_i$ for $i\le2$. This group generalizes the group $PSL(2,\mathbb{Z})$ to higher dimensions. The quotient space of $H^+$ with a principal subgroup $\Gamma_p[N]$ is indeed a manifold for $N \ge 2$, because $\Gamma_p[N]$ is torsion-free whenever $N \ge 2$. $\Gamma_p[N]$ consists of those matrices $\left(\begin{array}{cc} a & b \\ c & d \end{array}\right)$ from $\Gamma_p$ where the entries satisfy the arithmetic conditions $a-1,b,c,d-1 \equiv 0 $ mod $\Z + \Z e_1 + \Z e_2 + \Z e_3$. For more profound details on these groups and properties we refer the reader to \cite{BCKR}.  

\par\medskip\par

In order to generalize the representation formulas and the results that we obtained in the previous sections for the instationary Navier-Stokes system to the context of analogous instationary boundary value problems on conformally manifolds we only need to introduce the properly adapted analogues of the parabolic Dirac operator, and the other hypercomplex integral operators on these manifolds. So, the main goal consists of constructing the kernel functions explicitly. 
From the geometric point of view one is particularly interested in those conformally flat manifolds that have a spin structure, that means that one can construct at least one spinor bundle over such a manifold. These are called conformally spin manifolds. Often one can construct more than one spinor bundle over a spin manifold which leads to the consideration of spinor sections, in our case quaternionic spinor sections. For the geometric background we refer to \cite{LM}.  

\par\medskip\par

In this paper, we  restrict ourselves to explain the method at the simplest non-trivial example dealing with conformally flat spin $3$-tori with inequivalent spinor bundles. After this it becomes clear how to carry over our results to other examples of conformally flat (spin) manifolds that are constructed by factoring out a simply connected domain by a discrete Kleinian group, such as those mentioned above. 

\par\medskip\par
  
To start, let $\Omega:= \Z e_1 + \Z e_2 + \Z e_3$ be the standard lattice in $\mathbb{R}^3$. Then the topological quotient space $\R^3/\Omega$ is a $3$-dimensional conformally flat torus denoted by $T_3$, over which one can construct a number of conformally inequivalent spinor bundles over $T_3$.

We recall that in general different spin structures on a spin manifold $M$ are detected by the number of distinct homomorphisms from the fundamental group $\Pi_{1}(M)$ to the group ${\Z}_{2}=\{0,1\}$. In this case we have that $\Pi_{1}(T_3)={\Z}^{3}$. There are two homomorphisms of ${\Z}$ to ${\Z}_{2}$. The first one is $\theta_{1}:{\Z}\rightarrow {\Z}_{2}:\theta_{1}(n)=0$ mod $2$ while the second one is the homomorphism $\theta_{2}:{\Z}\rightarrow {\Z}_{2}:\theta_{2}(n)=1$ mod $2$. Consequently there are $2^{3}$ distinct spin structures on $T_3$. $T_{3}$ is a simple example of a Bieberbach manifold.  
\par\medskip\par
We shall now give an explicit construction for some of these spinor bundles over $T_3$. All the others are constructed similarly. First let $l$ be an integer in the set $\{1,2,3\}$, and consider the sublattice ${\Z}^{l}={\Z}e_{1}+\ldots+{\Z}e_{l}$  where$(0 \le l \le 3)$. In the case $l=0$ we simply have ${\Z}^{0}:=\emptyset$. 
There is also the remainder lattice  ${\Z}^{3-l}={\Z}e_{l+1}+\ldots+{\Z}e_{3}$. In this case ${\Z}^{3}=\{\underline{m}+\underline{n}:\underline{m}\in {\Z}^{l}$ and $\underline{n}\in {\Z}^{3-l}\}$.  Suppose now that $\underline{m}=m_{1}e_{1}+\ldots+m_{l}e_{l}$. Let us now make the identification $({\bf x},X)$ with $({\bf x}+\underline{m}+\underline{n},(-1)^{m_{1}+\ldots+m_{l}}X)$ where ${\bf x}\in \R^{3}$ and $X\in \mathbb{H}$. This identification gives rise to a quaternionic spinor bundle $E^{(l)}$ over $T_3$.

\par\medskip\par

Notice that $\R^3$ is the universal covering space of $T_3$. Consequently, there exists a well-defined
projection map $p: \R^3 \to T_3$. As explained for example in \cite{KraRyan2} every $3$-fold periodic resp. anti-periodic open set $U \subset \R^3$
and every $3$-fold periodic resp. anti-periodic section $f: U' \to E^{(l)}$, satisfying $f({\bf x}) = (-1)^{m_1+\cdots+m_l}({\bf x} + \omega)$ for all $\omega \in \Z^{l}\oplus \Z^{3-l}$, descends to a well-defined open set $U'=p(U) \subset T_3$ (associated with the chosen spinor bundle) and a well-defined spinor section $f':=p(f): U' \subset T_3 \to E^{(l)} \subset \mathbb{H}$, respectively. The projection map $p: \R^3 \to T_3$ induces well-defined toroidal modified parabolic Dirac operators on $T_3 \times \R^+$ by $p(D_{{\bf x},t,k}^{\pm}) =: {\cal {D}}_{{\bf x},t,k}^{\pm}$ acting on spinor sections of $T_3 \times \R^+$. Sections defined on open sets $U$ of $T_3 \times \R^+$ are called toroidal $k$-left parabolic monogenic if ${\cal {D}}_{{\bf x},t,k}^{\pm} s = 0$ holds in $U$. By ${\tilde{D}}:=p({\bf D})$ we denote the projection of the time independent Euclidean Dirac operator to the torus $T_3$. 

\par\medskip\par

The projections of the $3$-fold (anti-)periodization of the function $E({\bf x},t;k)$ denoted by 
$$ 
{\cal{E}}({\bf x},t;k) := \sum\limits_{\omega
\in \Z^3\oplus \Z^{3-l}}(-1)^{m_1+\cdots+m_l} E({\bf x}+\omega,t;k)
$$ 
provides us with the fundamental section to the toroidal parabolic modified Dirac operator 
${\cal {D}}_{{\bf x},t,k}^{\pm}$ acting on the corresponding spinor bundle of the torus $T_3$. From the function theoretical point of view the function ${\cal{E}}({\bf x},t;k)$ can be regarded as the canonical generalization of the classical elliptic Weierstra{\ss} $\wp$-function to the context of the modified Dirac operator $D_{{\bf x},t,k}^{+}$ in three dimensions. 

\par\medskip\par

To show that this expression is well-defined we have to prove the convergence of the series. So, the main task is show   
\begin{theorem} The series 
$$ 
{\cal{E}}({\bf x},t;k) = \sum\limits_{\omega \in \Z^3\oplus \Z^{3-l}}(-1)^{m_1+\cdots+m_l} E({\bf x}+\omega,t;k)
$$ 
converges uniformally on any compact subset of $\R^3 \times \R^+$. 
\end{theorem}
{\bf Proof}: We decompose the total lattice $\mathbb{Z}^3$ into the 
the following union of lattice points $\Omega = \bigcup_{m=0}^{+\infty} \Omega_m$ where
$$\Omega_m := \{ \omega \in \Z^3 \mid |\omega|_{max} = m\}.$$ 
We further consider the following subsets of this lattice 
$$
L_m := \{
\omega \in \Z^3 \mid |\omega|_{max} \le m\}.
$$
By a direct counting argument one observes the set $L_m$ contains exactly $(2m+1)^3$ many points. Hence, the
cardinality of $\Omega_m$ is $\sharp \Omega_m = (2m+1)^3 -
(2m-1)^3$. The Euclidean distance between the set $\Omega_{m+1}$
and $\Omega_m$ has the value $d_m := dist_2(\Omega_{m+1},\Omega_m)
= 1$.

\par\medskip\par

To show the normal convergence of the series, let us
consider an arbitrary compact subset ${\cal{K}} \subset \R^3$. Let $t > 0$ be an arbitrary but fixed value. 
Then there exists a positive real $r \in \R$ such that all ${\bf x} \in
{\cal{K}}$ satisfy $|{\bf x}|_{max} \le |{\bf x}|_2 < r$. 
Suppose now that ${\bf x}$ is a point of ${\cal{K}}$. To show the
normal convergence of the series we may leave out without loss of
generality a finite set of lattice points. So, we retrict ourselves to extend only the summation over those lattice points that
satisfy $|\omega|_{max} \ge [r]+1$.
In view of $$|{\bf x} + \omega|_2 \ge |\omega|_2 - |{\bf x}|_2 \ge
|\omega|_{max}-|{\bf x}|_2 = m - |{\bf x}|_2 \ge m - r$$ we
obtain
\begin{eqnarray*}
& & \sum\limits_{m=[r]+1}^{+\infty} \sum\limits_{\omega \in \Omega_m} |E({\bf x},t;k)({\bf x}+\omega)|_2\\
& \le & \frac{k}{(2 \sqrt{\pi t})^3} \sum\limits_{m=[r]+1}^{+\infty} \sum\limits_{\omega \in \Omega_m} \exp(-k|{\bf x}+\omega|_2/4t)\Big(\frac{k}{2t} |{\bf x} + \omega|_2 + \mathfrak{f}(\frac{3}{2t} 
+ \frac{k|{\bf x}+\omega|_2^2}{4t^2})+k \mathfrak{f}^{\dagger}\Big)\\
& \le & \frac{k}{(2 \sqrt{\pi t})^3} \sum\limits_{m=[r]+1}^{+\infty}\Big( [(2m+1)^3 - (2m-1)^3] \big(
\frac{k(r+m)}{2t} +\mathfrak{f}(\frac{3}{2t} + \frac{k(r+m)^2}{4t^2}) + k \mathfrak{f}^{\dagger}\big)\\
& & \quad\quad\quad\quad\quad\quad\quad\quad \times \;\; \exp(\frac{-k(m-r)^2}{4t})\Big),
\end{eqnarray*}
because $m - r \ge [r]+1-r > 0$. This sum clearly is absolutely uniformly convergent because of the decreasing exponent (remember $k > 0$) which dominates the polynomial expressions in $m$.  
Hence, the series 
$${\cal{E}}({\bf x},t;k) :=
\sum\limits_{\omega \in \Z^l \oplus \Z^{3-l}} (-1)^{m_1+\cdots+m_l} E({\bf x}+\omega,t;k),$$
which can be rewritten as
$$
{\cal{E}}({\bf x},t;k) :=
\sum\limits_{m=0}^{+\infty}\sum\limits_{\omega \in \Omega_m} (-1)^{m_1+\cdots+m_l}
E({\bf x}+\omega,t;k),
$$
converges normally on $\R^3 \times \R^+$. Since
$E({\bf x}+\omega,t;k)$ belongs to Ker $D_{{\bf x},t,k}^{+}$ in each
$({\bf x},t) \in \R^3 \times \R^+$  the series ${\cal{E}}({\bf x},t;k)$ satisfies
$D_{{\bf x},t,k}^{+} {\cal{E}}({\bf x},t;k)  = 0$ in 
each ${\bf x} \in \R^3 \times \R^+$. \hfill $\blacksquare$
\par\medskip\par
Obviously, by a direct rearrangement argument, one obtains that 
$$
{\cal{E}}({\bf x},t;k) =(-1)^{m_1+\cdots+m_l}  {\cal{E}}({\bf x}+\omega,t;k)\;\;\;\forall \omega \in \Omega
$$
which shows that the projection of this kernel correctly descends to a section with values in the spinor bundle $E^{(l)}$. The projection $p({\cal{E}}({\bf x},t;k))$ denoted by $\tilde{{\cal{E}}}({\bf x},t;k)$ is the fundamental section of the toroidal modified parabolic Dirac operator ${\tilde{D}}_{{\bf x},t,k}^{+}$. For a time-varying Lipschitz domain $G \subset T_3 \times \R^+$ with a strongly Lipschitz boundary $\Gamma$ we can now similarly introduce the Teodorescu and Cauchy transform for toroidal $k$-monogenic parabolic quaternionic spinor valued sections by 
\begin{eqnarray*}
\tilde{T}_G u({\bf y},t_0) &=& \int_G \tilde{{\cal{E}}}({\bf x}-{\bf y},t-t_0;k) u({\bf x},t) dV dt\\
\tilde{F}_{\Gamma} u({\bf y},t_0) &=& \int_{\Gamma} \tilde{{\cal{E}}}({\bf x}-{\bf y},t-t_0;k) d\sigma_{{\bf x},t} u({\bf x},t).
\end{eqnarray*}
Next, the associated Bergman projection can be introduced by 
$$
\tilde{{\bf P}} = \tilde{F}_{\Gamma}(tr_{\Gamma} \tilde{T}_G \tilde{F}_{\Gamma})^{-1} tr_{\Gamma}
\tilde{T}_G.
$$
and $\tilde{{\bf Q}} := \tilde{{\bf I}} - \tilde{{\bf P}}$. 
\par\medskip\par
Adapting from \cite{KraRyan2,CK2009} we obtain a direct analogy of Theorem~1, Lemma~1 and Theorem~2 on these conformally flat $3$-tori using these toroidal versions $\tilde{T}_G, \tilde{F}_{\Gamma}$ and $\tilde{{\bf P}}$ of operators introduced in Section~2. 
Suppose next that we have to solve a Navier-Stokes problem of the form (1)-(3) within a Lipschitz domain $G \subset T_3 \times \R^+$ with values in the spinor bundle $E^{(l)} \times \R^+$. Then we can compute its solutions by simply applying the following adapted iterative algorithm 
\begin{eqnarray*}
{\bf u}_n &=& \tilde{T}_G \tilde{{\bf Q}} \tilde{T}_G \Big[{\bf f} -\Re({\bf u}_{n-1}{\tilde {D}}){\bf u}_{n-1}\Big] - \tilde{T}_G \tilde{{\bf Q}}\tilde{T}_G {\tilde{D}} p_n\\
\Re(\tilde{{\bf Q}}\tilde{T}_G {\tilde{D}} p_n)&=& \Re\Big[\tilde{{\bf Q}} \tilde{T}_G 
{\bf f} - \Re({\bf u}_{n-1}{\tilde{D}}){\bf u}_{n-1}\Big]
\end{eqnarray*}
In the same flavor one obtains a direct analogy of Theorem~3 and Theorem~4 in this context. 
\par\medskip\par
Now it becomes clear how this approach even carries over to more general conformally flat spin manifolds that arise by factoring out a simply connected domain $U$ by a discrete Kleinian group $\Gamma$. The Cauchy-kernel is constructed by the projection of the $\Gamma$-periodization (involving eventually automorphy factors like in \cite{BCKR}) of the fundamental solution $E({\bf x};t;k)$. With this fundamental solution we construct the corresponding integral operators on the manifold. 

In terms of these integral operators we can express the solutions of the corresponding Navier-Stokes boundary value problem on these manifolds, simply by replacing the usual hypercomplex integral operators by its adequate ``periodic''  analogies on the manifold. This again underlines the very universal character of our approach to treat the Navier-Stokes equations but also many other complicated elliptic, parabolic, hypoelliptic and hyperbolic PDE systems with the quaternionic operator calculus using Dirac operators. 

Furthermore, the representation formulas and results also carry directly over to the $n$-dimensional case in which one simply replaces the corresponding quaternionic operators by Clifford algebra valued operators, such as suggested in \cite{CK2,CK2009}. 

\section{Acknowledgments}

The work of the third  author is supported by the project \textit{Neue funktionentheoretische Methoden f\"ur instation\"are PDE}, funded by Programme for Cooperation in Science between Portugal and Germany, DAAD-PPP Deutschland-Portugal, Ref: 57340281. The work of the first and second authors is supported via the project ``New Function Theoretical Methods in Computational Electrodynamics'' approved under the agreement A\c c\~oes Integradas Luso-Alem\~as DAAD-CRUP,  ref. A-15/17, and by Portuguese funds through the CIDMA - Center for Research and Development in Mathematics and Applications, and the Portuguese Foundation for Science and Technology (``FCT--Funda\c{c}\~ao para a Ci\^encia e a Tecnologia''),  within project UID/MAT/0416/2013.

\end{document}